\documentclass[12pt]{amsart}
\usepackage{amsthm, amssymb, url, enumerate, xypic, nicefrac}
\usepackage[all,cmtip]{xy}

\makeatletter
\@namedef{subjclassname@2010}{%
  \textup{2010} Mathematics Subject Classification}
\makeatother

\makeatletter
\let\@@pmod\pmod
\DeclareRobustCommand{\pmod}{\@ifstar\@pmods\@@pmod}
\def\@pmods#1{\mkern4mu({\operator@font mod}\mkern 6mu#1)}
\makeatother

\newtheorem{thm}{Theorem}[section]
\newtheorem{cor}[thm]{Corollary}
\newtheorem{lem}[thm]{Lemma}
\newtheorem{prop}[thm]{Proposition}

\theoremstyle{definition}
\newtheorem{defin}[thm]{Definition}

\newtheorem*{xrem}{Remark}

\DeclareMathOperator\Gal{Gal}
\DeclareMathOperator{\Ind}{Ind}
\def\cO{\mathcal{O}}
\def\cA{\mathcal{A}}
\def\fp{\mathfrak{p}}
\def\fa{\mathfrak{a}}
\def\fb{\mathfrak{b}}
\def\A{\mathbb{A}}

\def\Z{\mathbb{Z}}
\def\G{\mathbb{G}}
\def\N{\mathbf{N}}
\def\Q{\mathbb{Q}}
\def\R{\mathbb{R}}
\def\C{\mathbb{C}}

\def\GL{\operatorname{GL}}
\def\SL{\text{SL}}
\def\Res{\text{Res}}
\def\Hom{\text{Hom}}

\newcommand{\res}{\mathop{\rm res}}
\newcommand{\real}{\mathop{\rm Re}}
\newcommand{\eqdef}{\mathop{=}^{\rm def}}
\newcommand{\eps}{\varepsilon}
\newcommand{\es}[1]{\begin{equation}\begin{split}#1\end{split}\end{equation}}
\newcommand{\est}[1]{\begin{equation*}\begin{split}#1\end{split}\end{equation*}}
\newcommand{\legen}[2]{\left(\frac{#1}{#2}\right)}
\newcommand{\bijar}[1][]{%
 \ar[#1]
 \ar@<0.7ex>@{}[#1]|-*=0[@]{\sim}}

\begin{document}
\title{Counting automorphic forms on norm one tori}

\baselineskip=17pt

\author[E.~Hunter Brooks]{Ernest Hunter Brooks} 
\address{Section de Math\'ematiques\\
	Ecole Polytechnique F\'ed\'erale de Lausanne\\
	B\^atiment MA, Station 8 \\
	1015 Lausanne \\
	Switzerland}
\email{ernest.brooks@epfl.ch}

\author{Ian Petrow} 
\address{Section de Math\'ematiques\\
	Ecole Polytechnique F\'ed\'erale de Lausanne\\
	B\^atiment MA, Station 8 \\
	1015 Lausanne \\
	Switzerland}
\curraddr{ETH Z\"urich - Departement Mathematik \\
HG G 66.4 \\
R\"amistrasse 101 \\
8092 Z\"urich \\
Switzerland}
\email{ian.petrow@math.ethz.ch}

\date{}

\begin{abstract}
We give an asymptotic formula for the number of automorphic forms on the non-split norm one torus $T$ associated with an imaginary quadratic extension of $\mathbb{Q}$, ordered by analytic conductor.  \end{abstract}



\maketitle

\section{Introduction}

In the study of automorphic forms, an interesting question is whether they are more abundant on some groups than on others. To study this phenomenon precisely, it is necessary to count automorphic forms, which requires ordering them by some positive real invariant with the property that the number of automorphic forms with invariant below a fixed bound is finite. To this end, Sarnak, Shin, and Templier \cite{SarnakShinTemplier} have proposed using the \textit{analytic conductor}, motivated by an analogy between the set of automorphic forms of bounded analytic conductor and the set of points of bounded \textit{height} on a variety over a global field.

Accordingly, for a connected reductive group $G$, write $\cA(G)$ for the set of all automorphic forms on $G$ and consider \est{f_G(X) = \#\{ \pi \in \cA(G)  \text{ with analytic conductor}\leq X \}.}  In recent work, Brumley and Mili\'cevi\'c \cite{BrumleyMilicevic} have computed an asymptotic formula for $f_G(x)$ when $G = \GL_1$ or $G = \GL_2$ over an arbitrary number field, confirming a conjecture of Michel and Venkatesh for these cases \cite[p.~204]{MichelVenkateshGL2}. In this paper we give a formula for the case where $G=T$ is the non-split one-dimensional torus over $\Q$ associated with an imaginary quadratic field $K$.

The analytic conductor, first introduced by Iwaniec and Sarnak (see \cite[Eq.~31]{IS2000GAFA}), is a product of local terms which can be read off of the completed global $L$-function attached to an automorphic form $\pi  \in \cA(\GL_n)$.  To define analytic conductors (and $L$-functions) for automorphic forms on more general reductive groups $G$, one must choose a finite-dimensional complex representation of the $L$-group of $G$. For an arbitrary group $G$ there is no privileged choice of representation, however for $G = T$ there is a natural such choice, as explained in  \textsection 2. The analytic conductor must also be normalized at archimedean places $v \mid \infty$.  This paper takes the convention that the local analytic conductor at $v\mid \infty$ is $C(s,\psi_v)=(|s+\kappa_v|+1)^{[K_v:\R]},$ where $[K_v:\R] =1 \text{ or }2$.

Write $h_K,w_K$ and $d_K$ for the class number, number of units and discriminant of $K$. Write also $C(\psi)=C(0,\psi)$.

\begin{thm}\label{MT}  
For $X\geq3|d_K|$, 
\est{ f_T(X) = \frac{h_K}{w_K}\frac{2\zeta(2)-1}{\zeta_K(2)} \frac{X}{|d_K|}+ O( h_K \left(\frac{X}{|d_K|}\right)^{2/3}(\log X/|d_K|)^3) ,} and
\est{\sum_{\substack{\psi \in \cA(T) \\ \psi(K_\infty^\times ) = 1 \\ C(\psi)\leq X}} 1= \frac{h_K}{w_K} \frac{1}{\zeta_K(2)} \frac{X}{|d_K|}+ O( h_K \left(\frac{X}{|d_K|}\right)^{2/3}(\log X/|d_K|)^2) .} The implied constants are absolute, i.e., they do not depend on $K$. \end{thm}  Theorem \ref{MT} follows, after a generating series computation, from the following more precise result.  For $p$ a rational prime, write \est{\phi^*(p^j) = \begin{cases} 
1 & \text{ if } j=0 \\ 
0 & \text{ if } j \text{ is odd} \\
p-2 & \text{ if } j=2 \text{ and } p \text{ splits} \\ 
p^{\frac{j}{2}-2}(p-1)^2 & \text{ if } j\geq 4, \text{ is even, and } p \text{ splits} \\ 
p & \text{ if } j=2 \text{ and } p \text{ is inert} \\ 
p^{\frac{j}{2}-2}(p^2-1) & \text{ if } j\geq 4, \text{ is even, and } p \text{ is inert} \\  
p^{\frac{j}{2}-1}(p-1) & \text{ if } j\geq 2, \text{ is even, and } p \text{ is ramified},  
\end{cases}} and extend the function $\phi^*$ to all $n \in \mathbb{N}$ multiplicatively.  
\begin{thm}\label{MT2}
\est{\sum_{\substack{\psi \in \cA(T) \\ \psi(K_\infty^\times ) = 1\\ C(\psi)=n|d_K| }} 1 = \begin{cases} h_K \phi^*(n) & \text{ if } K \neq \Q(i),\Q(\zeta_3) \\ h_K\left(\frac{1}{2}\phi^*(n)+ \frac{1}{2} \mu(\sqrt{n})\right) & \text{ if } K=\Q(i) \\  h_K\left(\frac{1}{3}\phi^*(n)+ \frac{2}{3} \mu(\sqrt{n})\right) & \text{ if } K=\Q(\zeta_3),\end{cases}} where $\mu(m)$ is the M\"obius function and $\mu(\sqrt{n})=0$ if $n$ is not a perfect square.
\end{thm}
Remarks: \begin{enumerate}[\upshape (i)]
\item The first statement in Theorem \ref{MT} above depends on the aforementioned choice of definition of the analytic conductor at the archimedean place, whereas the second statement does not.
\item Recall the classical analytic class number formula for $K$ an imaginary quadratic field of discriminant $d<0$: $h_K/w_K = (2\pi)^{-1}\sqrt{|d|}L(1,\chi_{d}),$ where $\chi_d(n) = \legen{d}{n}$ is given by the Kronecker symbol.  The leading constants in Theorem \ref{MT} then admit the expected factorizations over the places of $\Q$.
\item Theorem \ref{MT2} counts only those $\psi \in  \cA(T)$ with trivial infinity type.  The proof in  \textsection \ref{GeneratingSeries} also gives a count for any other specified infinity type, with a similar formula.  
\item The special value $\zeta_K(2)$ is not critical in the sense of Deligne and thus is not expected to admit a closed-form description in terms of familiar constants. A theorem of Borel relates it to the volume of $\SL(2, \C)/\SL(2, \cO_K)$ for Tamagawa's normalization of the Haar measure and to the third $K$-group of $K$ (see \cite{BorelKtheory}).
\item It is interesting to compare these results to the case of a split torus over $\Q$, that is, the case $\GL_{1}/\Q$.  The natural analogue of Theorem \ref{MT} is to count $\psi \in \cA(\GL_{1})$ according to the conductor defined via the representation $z \oplus z^{-1}$ of the $L$-group $\C^\times$ of $\GL_1$.  The automorphic form associated to $\psi$ by functoriality is the Eisenstein series $\psi \oplus \psi^{-1}$, whose conductor is $C(\psi)^2$, where $C(\psi)$ is the classical conductor of $\psi$.  Counting via this conductor, a classical computation gives \est{ \sum_{\substack{\psi \in \cA(\GL_{1}/\Q) \\ \psi(\R^\times)=1 \\ C(\psi \oplus \psi^{-1}) \leq X}} 1 \sim \frac{1}{2\zeta(2)^2}X .}  \item The implied constant in the error term in Theorem \ref{MT} is effectively computable.  The error term in Theorem \ref{MT} could be improved slightly, but at the cost of introducing an ineffective dependence on $K$ in the implied constant owing to the fact that one cannot rule out the existence of Siegel zeros for the Dirichlet $L$-functions $L(s, \chi_d)$.
\item A similar result to Theorem \ref{MT} is expected to hold for the non-split one-dimensional torus over $\Q$ splitting over a \emph{real} quadratic field. In this case, the space $\cA(T)$ is no longer discrete and one needs to be more careful in setting up the problem. 
\end{enumerate}

For the purpose of comparison to Theorem \ref{MT}, we recall below the results of Brumley and Mili\'cevi\'c on $\GL_1$ and $\GL_2$.  The counting problem for inner forms of $\GL_2$ has also been studied by Lesesvre.

Let $F$ be a number field and let $\mathfrak{F}$ be the set of unitary cuspidal automorphic representations $\pi$ of $\GL_n(\A_F)$, such that the central character of $\pi$ is trivial on $\R_{>0}$, diagonally embedded in $\A_F^\times$ as $t \mapsto \prod_{v \mid \infty} t^{1/[F:\Q]}$.  Brumley and Mili\'cevi\'c choose this normalization to eliminate continuous families of automorphic forms. For any finite place $v$ of $F$ let $$\Delta_v(s) = \zeta_v(s) \zeta_v(s+1) \cdots \zeta_v(s+n-1).$$ Let $\Delta_F(s) = \prod_v \Delta_v(s)$ and write $\Delta^*_F(1)$ for its residue at $s=1$. Write also $\zeta_F^*(1)$ for the residue of $\zeta_F(s)$ at $s=1$.  Let $$\lVert \widehat{\nu}(\mathfrak{F}) \rVert = \frac{\zeta_F^*(1)}{\zeta_F(n+1)^{n+1}}\lVert\widehat{\nu}_\infty \rVert,$$ where $$\lVert \widehat{\nu}_\infty\rVert = \int_{\Pi(\GL_{n}(F_\infty)^1)} q(\pi_\infty)^{-n-1} \,d\widehat{\mu}_{\infty}^{\rm pl}(\pi_\infty).$$ Here, $\Pi(\GL_{n}(F_\infty))$ is the local unitary dual at infinity, the superscript $1$ means the subset where the above mentioned normalization on the central character holds, and $\widehat{\mu}_{\infty}^{\rm pl}$ is the Plancherel measure on $\Pi(\GL_{n}(F_\infty))$, normalized so that the Plancherel formula holds relative to the canonical measure $\Delta_\infty(1)\det(g)^{-n}(dg_{11} \wedge \cdots \wedge dg_{nn})$. 

Let $F$ be fixed and $n=1$ or $n=2$.  Brumley and Mili\'cevi\'c prove that as $X \to \infty$ we have $$|\{ \pi \in \mathfrak{F} : Q(\pi) \leq X\}| \sim \frac{1}{n+1} |d_F|^{n^2/2} \Delta_F^*(1)\lVert \widehat{\nu}(\mathfrak{F})\rVert X^{n+1},$$ where $Q(\pi)$ is the global analytic conductor of $\pi$.

To prove Theorem \ref{MT}, in  \textsection \ref{AC} we define the analytic conductor for forms on $T$ by calculating their $L$-functions; the main difficulty here is the calculation of Euler factors at ramified primes, which is done automorphically and Galois-theoretically,  with part of the latter calculation taking place in an appendix.  In  \textsection \ref{setup}, the $\psi \in \cA(T)$ and their analytic conductors  are related to more readily countable objects.  Finally,  \textsection \ref{GeneratingSeries} contains counting arguments to arrive at Theorem \ref{MT2}, from which Theorem \ref{MT} follows from a generating series computation.

\section{$L$-functions and the Analytic Conductor}\label{AC}

\subsection{$L$-functions}\label{Lfcns}

To give a definition of the analytic conductor for characters $\psi \in \cA(T)$, we first study their $L$-functions. If one were only interested in counting on $T$ and not on more general groups, one could define an $L$-function for $\psi$ in an essentially ad hoc fashion (by taking Proposition \ref{sharpmap} below as a definition). However, in order to meaningfully compare the growth rates in Theorem \ref{MT} to those for other groups, it is necessary to take an approach which generalizes.

Given an arbitrary reductive group $G$ over $\Q$, an automorphic form $\pi$ on $G$, and a representation $\rho: {}^L G \to \GL_n(\C)$ of the Langlands dual group of $G$, one has a recipe for defining Euler factors $E_p(s, \pi, \rho)$ at all but finitely many places $p$ (see \cite[p. 10]{ArthurGelbart}, or \cite[p. 49]{BorelAutomorphicLfcns}, or read  \textsection \ref{goodplaces} for the application of this recipe to the specific case $G = T$). The analytic conductor depends crucially on the data at \emph{all} places $p$, and so one needs a definition of the correct Euler factors at other places.

There are (at least) two possible approaches in general, both of which depend on significant auxiliary results for $G$, and which (in this generality) only conjecturally yield the same answer:
\begin{itemize}
\item If one knows the global Langlands correspondence for $G$ (see \cite[Chap. III]{BorelAutomorphicLfcns} or  \textsection \ref{artinL} of this paper for definitions), then there is a Langlands parameter $\alpha_\pi: W_\Q \to {}^L G$ attached to $\pi$, where $W_\Q$ denotes the Weil (or Weil-Deligne) group of $\Q$, and one may compose $\alpha_\pi$ with $\rho$ to define an Artin $L$-function, whose Euler factors will be defined at all places and which will agree with $E_p(s, \pi, \rho)$ at the set of places where the latter are defined. This is the definition given in \cite[\textsection 12.1]{BorelAutomorphicLfcns}, for example.
\item The representation $\rho$ induces a unique $L$-homomorphism ${}^L G \to {}^L \GL_n$. If one knows the Langlands functoriality conjecture (see \cite[\textsection 17]{BorelAutomorphicLfcns} or  \textsection \ref{functoriality} of this paper) for this morphism, then there is an automorphic form $\rho_*(\pi)$ on $\GL_n$ for which there is a purely automorphic construction of Euler factors at all places, and these agree with the previously defined local factors where both are defined. By multiplicity one for $\GL_n$, the missing Euler factors are then uniquely determined from the ones that can be calculated directly. This is the approach taken in \cite{KowalskiFamilies}.
\end{itemize}

Either approach works for tori, because in this case the global Langlands correspondence is a theorem of Langlands (\cite{LanglandsAbelian}), and functoriality can be shown by a combination of simple arguments for tori and known results for two-dimensional representations of ${}^L \text{Res}_{K/\Q} \G_m$. We take the latter approach, but also calculate the relevant Artin $L$-functions in  \textsection \ref{artinL} to show that the definitions agree.

It should be noted that the analytic conductor of an automorphic form can, given a choice of $\rho$ as above, be defined as a product of local factors without explicit mention of the theory of $L$-functions, as in \cite[1.9]{KowalskiFamilies}.

\subsection{The torus $T$ and its $L$-group}\label{torusfacts}
The non-split torus $T$ is the kernel of the norm map $\text{Res}_{K/\Q} \G_m \to \G_m$.

If $K = \Q(\sqrt{d})$ and $p$ splits in $K/\Q$, then there is an $x \in \Q_p$ such that $x^2 = d$, and one has an identification of $\Q_p$ algebras $$K \otimes \Q_p = \Q_p \times \Q_p$$ given by $a + b\sqrt{d} \mapsto (a + b x, a - bx).$ Under this identification, the generator $\sigma$ of $\Gal(K/\Q)$ identifies with the involution of $\Q_p \times \Q_p$ given by $(\alpha, \beta) \mapsto (\beta, \alpha)$, so the norm map identifies with the product map $\Q_p^\times \times \Q_p^\times \to \Q_p^\times$ and $T(\Q_p)$ identifies with the set of anti-diagonal elements of $\Q_p^\times \times \Q_p^\times$. Its maximal compact subgroup corresponds under this identification to the set of anti-diagonal elements of $\Z_p^\times \times \Z_p^\times$.

If $p$ is non-split in $K/\Q$, then $T(\Q_p)$ is the set of norm one elements in the field $K \otimes \Q_p$, and in particular $T(\Q_p)$ is compact. 

The character lattice of $T$ is isomorphic to $\mathbb{Z}$, with the non-trivial element of $\Gal(K/\Q)$ acting by $-1$. It follows that ${}^LT= \C^\times \rtimes \Gal(K/\Q)$, where $\Gal(K/\Q)$ acts by inversion on $\C^\times$ (see \cite[\textsection 2]{BorelAutomorphicLfcns} for the general definition of $L$-groups; our definition is the variant in Remark 2.4.(2) of \textit{loc.~cit.}). A complete list of isomorphism classes of irreducible finite-dimensional complex algebraic representations of ${}^L T$ is given by:
\begin{itemize}
\item The trivial representation $\C$.
\item The sign representation $\C(-)$, the one-dimensional representation on which ${}^L T$ acts via ${}^L T \to \Gal(K/\Q)$.
\item For $n \geq 1$, the $2$-dimensional representation $V_n := \Ind_{\C^\times}^{{}^L T} \C(n)$, where $\C(n)$ is $\C$ with $z \in \C^\times$ acting by $z^n$. A model for $V_n$ is given by $z \rtimes 1 \mapsto \begin{pmatrix} z^n & 0 \\ 0 & z^{-n} \end{pmatrix}$, and $1 \rtimes \sigma \mapsto \begin{pmatrix} 0 & 1 \\ 1 & 0 \end{pmatrix}$.
\end{itemize}
We refer to $V_1$ as the standard representation (see  \textsection \ref{standardrep}).

\subsection{Euler factors at unramified places}\label{goodplaces}

For each prime $p$ unramified in $K/\Q$, the Satake correspondence gives a bijection between the unramified characters $\psi_p$ of $T(\Q_p)$ (for $T$, these are exactly the characters which are trivial on the maximal compact subgroup) and the conjugacy classes in ${}^LT$ whose projection onto the factor $\Gal(K/\Q)$ is equal to the Frobenius at $p$  (see \cite[\S 9.5]{BorelAutomorphicLfcns}). To write this correspondence explicitly for $T$, write $\psi \in \cA(T)$ as $\psi = \bigotimes_v \psi_v$.

The set $M_{\text{big}}=\{ z \rtimes \sigma : z \in \C^\times\} \subset {}^L T$ forms a single conjugacy class. The other conjugacy classes in ${}^LT$ are given by the sets $M_z = \{z \rtimes 1,z^{-1} \rtimes 1\}$.  

If $p$ is inert in $K/\Q$ and $\psi$ is unramified at $p$, then $\psi_p$ is necessarily trivial as $T(\Q_p)$ is compact, corresponding to the conjugacy class $M_\text{big}$.  

At split places $p$, an unramified character $\psi_p$ is determined by its image on $(p,p^{-1}) \in T(\Q_p)$, and the correspondence is given by \est{ \psi_p \leftrightarrow M_{\psi_p((p,p^{-1}))}.}

Given a representation $\rho:{}^LT \to \GL_m(\C)$ and a $\psi \in \cA(T)$, unramified at $p$, the Euler factor $E_p(s, \psi, \rho)$ at $p$ is defined to be the reciprocal of the (co-)characteristic polynomial of the image under $\rho$ of any element of the conjugacy class corresponding to $\psi_p$, evaluated at $p^{-s}$. These polynomials are recorded in the table below for the various irreducible representations of ${}^LT$.
\begin{center}
\begin{tabular}{ l | l | l }
Rep.~$\rho$& Char.~poly.~at split $p$ & Char.~poly.~at inert $p$  
\\ \hline \hline
$\C$ & $1 - x$ & $1 - x$   
\\ \hline
$\C(-)$ & $1 - x$ & $1 + x$  
\\ \hline
$V_n$ & $(1 - \psi_p(p, p^{-1})^n x)(1 - \psi_p(p, p^{-1})^{-n} x )$ & $1 -x^2$  
\\
\end{tabular}
\end{center}

Let $P$ be the set of places of $\Q$ given by the union of the ramified places of $K/\Q$, the ramified places of $\psi$, and the archimedean place $\infty$.  From the table, it is clear that, away from Euler factors at places in $P$, the various $L$-functions of $\psi$ (as $\rho$ varies) identify with classical $L$-functions. More precisely, given $\psi \in \cA(T)$, define a Hecke character $\psi^\sharp \in \cA(\GL_{1}/K)$ by the rule \est{\psi^\sharp(x) \eqdef \psi(x/x^\sigma).}  Then
\begin{center}
 \begin{tabular}{ l | l}
Rep.~$\rho$& $L_P(s,\psi,\rho)$ \\ \hline \hline 
$\C$ & $\zeta_P(s)$ \\ \hline
$\C(-)$ & $L_P(s,\chi_K)$ \\ \hline 
$V_n$ & $L_P(s,(\psi^\sharp)^n)$ \\
\end{tabular}
\end{center}
where the subscript $P$ means that in the Euler product defining each of the above $L$-functions the factors at primes dividing any prime in $P$ are omitted, $\zeta(s)$ denotes the Riemann zeta function, and $L(s,\chi_K)$ denotes the Dirichlet $L$-function of the quadratic character associated with the extension $K/\Q$. 

\subsection{The torus $S$ and functoriality}\label{functoriality}
The above table makes it easy to guess the missing Euler factors for each representation $\rho$ of ${}^L T$; namely, they coincide with the ones attached to the corresponding classical $L$-functions. As explained in \textsection\ref{Lfcns}, to make this rigorous, one must establish, given a choice of $\rho$, the functoriality conjecture for the corresponding $L$-homomorphism ${}^L T \to {}^L \GL_n$. (Here, ${}^L \GL_n$ denotes the direct product $\GL_n(\C) \times \Gal(K/\Q)$, an $L$-homomorphism is one commuting with the projections to $\Gal(K/\Q)$, and one sees easily that a homomorphism ${}^L T \to  \GL_n$ extends uniquely to an $L$-homomorphism  ${}^L T \to {}^L \GL_n$.)

Recall (see e.g. \cite[pp. 11--12]{ArthurGelbart}) that an $L$-homomorphism of $L$-groups ${}^L r: {}^L G \to {}^L H$ is said to satisfy the functoriality conjecture if there is a map $\sharp: \cA(G) \to \cA(H)$ such that for each $\pi \in \cA(G)$ and each representation $\rho$ of ${}^L H$, one has, for some finite set $P$ of places,
$$
L_P(s, \pi, \rho \circ {}^L r) = L_P(s, \pi^\sharp, \rho).
$$

Write $S = \Res_{K/\Q} \G_m$. We will show that all of the above representations factor through a universal inclusion ${}^L T \hookrightarrow {}^L S$, which satisfies the functoriality conjecture, and conclude by appealing to known functoriality results for a particular embedding ${}^L S \to {}^L \GL_2$ which induces the standard representation of ${}^L T$. 

We follow the same steps as above to compute the $L$-functions of forms $\chi \in \cA(S)$ and representations $\rho$ of ${}^L S$. Observe that ${}^LS = (\C^\times)^2\rtimes \Gal(K/\Q)$, with $\Gal(K/\Q)$ acting by interchanging the factors. The following is a complete list of its irreducible finite-dimensional complex algebraic representations:
\begin{itemize}
\item For $m \in \Z$, the one-dimensional representations $\C(m)$ on which $\widehat{S}$ acts via $(z_1z_2)^m$ and $\Gal(K/\Q)$ acts trivially.
\item For $m \in \Z$, the one-dimensional representations $\C(m)(-)$ on which $\widehat{S}$ acts via $(z_1z_2)^m$ and $\Gal(K/\Q)$ acts via negation.
\item For $m > n \in \Z$, the representation $V_{m, n} := \Ind_{\widehat{S}}^{{}^L S} \C(m, n)$, where $\C(m, n)$ is $\C$ with $(z_1, z_2)$ acting by $z_1^mz_2^n$. Explicitly, $(z_1, z_2)\mapsto \begin{pmatrix} z_1^mz_2^n & 0 \\ 0 & z_1^nz_2^m \end{pmatrix}$ and $\sigma \mapsto \begin{pmatrix} 0 & 1 \\ 1 & 0 \end{pmatrix}.$
\end{itemize}
Writing the two elements of $\Gal(K/\Q)$ as $1,\sigma$, the conjugacy classes within ${}^LS$ are given by 
\begin{itemize}
\item $M_z=\{(z_1,z_2) \rtimes \sigma : z_1z_2=z\}$
\item $M_{z_1,z_2}= \{(z_1,z_2) \rtimes 1, (z_2,z_1)\rtimes 1\}$.
\end{itemize}
To a prime $p$ unramified in $K/\Q$ and an unramified character $\chi_p$ of $S(\Q_p)$, one assigns the following conjugacy class of ${}^LS$:
\begin{itemize}
\item If $p=\fp\overline{\fp}$ splits then $\chi_p \leftrightarrow M_{\chi_p(\pi),\chi_p(\overline{\pi})}$, where $\pi,\overline{\pi}$ are local uniformizers of $\fp$, $\overline{\fp}$, respectively.
\item If $p$ is inert then $\chi_p \leftrightarrow M_{\chi_p(p)}$.
\end{itemize}
The tables below list the characteristic polynomials and partial $L$-functions for $S$:
\begin{center}
\begin{tabular}{ l | l | l}
Rep.~$\rho$ & Char.~poly.~at split $p = \pi\overline{\pi}$ & Char.~poly.~at inert $p$  \\ \hline \hline
$\C(m)$ & $1 - (\chi(\pi)\chi(\overline{\pi}))^m x$ & $1 - \chi(p)^m x$  \\
$\C(m)(-)$ & $1 - (\chi(\pi)\chi(\overline{\pi}))^m x$  & $1 + \chi(p)^m x$  \\
$V_{m, n}$ & $(1 - \chi(\pi)^m\chi(\overline{\pi})^nx)(1 - \chi(\pi)^n\chi(\overline{\pi})^mx)$ & $(1 - \chi(p)^{m+n}x^2)$  \\
\end{tabular}
\begin{tabular}{ l | l}
Rep.~$\rho$ &  $L_P(s,\chi,\rho)$ \\ \hline \hline
$\C(m)$ & $L_P(s, \chi^m|_{\A_\Q^\times})$ \\
$\C(m)(-)$ & $L_P(s, \chi_K \cdot \chi^m|_{\A_\Q^\times})$ \\
$V_{m, n}$  & $L_P(s, \chi^m \otimes {\chi^\sigma}^n)$  \\
\end{tabular}
\end{center}

Let  ${}^L r: {}^L T \to {}^L S$ be given by $z \rtimes \sigma^\epsilon \mapsto (z, z^{-1}) \rtimes \sigma^\epsilon$. Under this representation, $\C(m)$ restricts to the trivial representation of ${}^L T$, $\C(m)(-)$ to the sign representation, and $V_{m, n}$ to $V_{m-n}$.

Explicit comparison of the tables for $T$ and $S$ establishes functoriality for the morphism ${}^L r$, where the association $\sharp: \cA(T) \to \cA(S)$ is given by $\psi^\sharp(x) = \psi(x/x^\sigma)$ (note that $\psi^\sharp|_{\A_\Q^\times}$ is trivial). The morphism ${}^L r$ is \emph{universal} in the sense that every irreducible representation of ${}^L T$ arises as a restriction of an irreducible representation of ${}^L S$ along this embedding. Moreover, it is the unique $L$-homomorphism from ${}^L T$ to ${}^L S$ with this property.

\begin{xrem}
If one is not interested in the explicit $L$-functions, one may more easily prove functoriality for ${}^L r$ as follows: the map ${}^L r$ is induced by a map $r: S \to T$, given on $A$-points by $x \mapsto x/x^\sigma$. The correspondence $\sharp$ is pullback of characters along $r$, and functoriality follows formally from the functoriality of the Satake correspondence. The same argument generalizes to any morphism of tori.
\end{xrem}

To complete the proof that the morphism ${}^L T \to {}^L \GL_2$ given by the standard representation satisfies functoriality, we need to show that the morphism of $L$-groups induced by the representation $V_{1, 0}$ of $S$ satisfies the functoriality conjecture. This is the principle of automorphic induction, the $n=1$ case of \cite[Example 1b, p. 12]{ArthurGelbart}. It then follows that for $\rho$ the standard representation of ${}^LT$, the $L$-function $L(s, \psi, \rho)$ with all finite places is the usual Hecke $L$-function $L(s, \psi^\sharp)$.

\subsection{Comparison with Artin $L$-functions}\label{artinL}
For the purpose of generalizations to other groups where a proof of functoriality may not be available, and to verify that the two definitions of  \textsection \ref{Lfcns} agree, we quickly show that the $L$-factors computed at ramified places above agree with those of the Artin $L$-functions attached to $\psi \in \cA(T)$. For more detail on Weil groups and Artin $L$-functions, see \cite{tateNTB}. 

Recall the Weil group $W_\Q$, its quotient
$$
W_{K/\Q} = W_\Q/[W_K, W_K],
$$
and the exact sequence
\begin{equation}\label{weilgroup}
1 \to C_K \to W_{K/\Q} \to \Gal(K/\Q) \to 1,
\end{equation}
where $C_K$ is the id\`{e}le class group of $K$. The conjugation action of $\Gal(K/\Q)$ on $C_K$ induced by (\ref{weilgroup}) coincides with the usual action. For each $p$ there is a local Weil group $W_{\Q_p}$ and a restriction map, well-defined up to conjugacy, $W_{\Q} \to W_{\Q_p}$. The local Weil group sits in an exact sequence
$$
1 \to \mathcal{I}_p \to W_{\Q_p} \to \Z \to 0,
$$
where $\mathcal{I}_p$ is the inertia subgroup of the absolute Galois group of $\Q_p$ and any lift of $1 \in \Z$ is called a Frobenius element. 

A Langlands parameter for $T$ is a map $W_\Q \to {}^LT$ inducing the identity map on their common quotient $\Gal(K/\Q)$. Composition with the (set-theoretic) projection ${}^LT \to \widehat{T}$ gives a 1-cocycle $\xi \in H^1(W_\Q, \widehat{T})$. Say that two Langlands parameters are equivalent if the associated cocycles are everywhere locally equal in cohomology, i.e. if their difference restricts to a coboundary in $H^1(W_{\Q_p}, \widehat{T})$ for all $p$.

Given a Langlands parameter $\alpha$ and a representation $(\rho, V)$ of ${}^L T$, one has Artin $L$-functions $L_{Art}(s, \alpha, \rho)$. The Euler factor at $p$ is the reciprocal of the (co-)characteristic polynomial of a Frobenius element of $W_{\Q_p}$ acting via $\rho \circ \alpha$ on the invariants of $V$ for $\mathcal{I}_p$, evaluated at $p^{-s}$.

The Langlands correspondence, which is a theorem of Langlands in the case of tori, attaches a Langlands parameter (well-defined up to equivalence) $\alpha_\psi$ to any automorphic form $\psi$ on any torus (see \cite[Theorem 2.b]{LanglandsAbelian}). In the case of $T$, one has the following:

\begin{lem}\label{cocycle}
Choose a lift $\widetilde{\sigma}$ of the generator $\sigma$ of $\Gal(K/\Q)$ to $W_{K/\Q}$. For $\psi \in \cA(T)$, the Langlands parameter $\alpha_\psi: W_\Q \to {}^LT$ is given by the following formula:
\begin{itemize}
\item It is trivial on the kernel $[W_K, W_K]$ of the map $W_\Q \to W_{K/\Q}$.
\item For $x \in C_K \subset W_{K/\Q}$, one has 
$$\alpha_\psi(x) = \psi(x/x^\sigma).$$
\item For $x = \widetilde{\sigma} y \in  W_{K/\Q} \setminus C_K$, one has 
$$\alpha_\psi(x) = \psi(y^\sigma/y) \rtimes \sigma.$$
\end{itemize}
(A different choice of $\widetilde{\sigma}$ does not change the equivalence class of the Langlands parameter.)
\end{lem}
\begin{proof}
This follows from the explicit construction of \cite{LanglandsAbelian}. Because the calculations involve notation that will not be used elsewhere in the paper, we have relegated them to the appendix.
\end{proof}

\begin{cor}
For any character $\psi \in \cA(T)$ and $\rho$ the standard representation of ${}^LT$, the Artin $L$-function $L_{Art}(s, \alpha_\psi, \rho)$ coincides with the Hecke $L$-function $L(s, \psi^\sharp)$.
\end{cor}

\begin{proof}

Let $\psi$ be a character of $T(\A)/T(\Q)$, and write $\rho_\psi$ for $\rho \circ \alpha_\psi$, where $\rho$ is the standard representation.

By the formalism of Artin $L$-functions, it suffices to give an isomorphism of representations
$$
\rho_\psi = \Ind_{W_K}^{W_\Q} \psi^\sharp.
$$
where $\psi^\sharp$ is viewed as a representation of $W_K$ acting via its abelianization $C_K$. Both representations being trivial on $[W_K, W_K]$, it then suffices to give an isomorphism of representations
$$
\overline{\rho_\psi} = \Ind_{C_K}^{W_{K/\Q}} \psi^\sharp,
$$
where $\overline{\rho_\psi}: W_{K/\Q} \to \text{Aut}(V_1)$ is the representation whose lift to $W_K$ is $\rho_\psi$.
Consider the model for $V_1$ of \textsection\ref{torusfacts}. By Lemma \ref{cocycle}, under $\overline{\rho_\psi}$, for $x \in C_K$, one has
$$
x \mapsto \begin{pmatrix}\psi^\sharp(x) & 0 \\ 0 & \psi^\sharp(x)^{-1} \end{pmatrix},
$$
and for $x \in W_{K/F} \setminus C_K$, one has
$$
x \mapsto \begin{pmatrix} \psi^\sharp(\widetilde{\sigma}^{-1}x)^{-1} & 0 \\ 0 & \psi^\sharp(\widetilde{\sigma}^{-1}x)  \end{pmatrix} \begin{pmatrix} 0 & 1 \\ 1 & 0 \end{pmatrix}= \begin{pmatrix} 0 & \psi^\sharp(\widetilde{\sigma}^{-1}x)^{-1} \\ \psi^\sharp(\widetilde{\sigma}^{-1}x) & 0 \end{pmatrix}.
$$
This is exactly the usual realization of the induced representation of $\psi^\sharp$.
\end{proof}

\subsection{Remarks on the choice of representation $V_1$}\label{standardrep}
For the purpose of generalization to other groups, we list the reasons to choose the particular representation $V_1$ of ${}^L T$:
\begin{itemize}
\item The representation $V_{1, 0}$ of ${}^L S$ is the representation that occurs in the automorphic induction theorem which allows one to view Hecke characters of $K$ as $\GL_2$ $L$-functions, and its restriction under the universal embedding ${}^L T \hookrightarrow {}^L S$ is $V_1$.

\item For any non-exceptional group, Gan, Gross, and Prasad give a conjectural standard representation of the $L$-group. The standard representation of this paper coincides with theirs (\cite[pp.25-26]{GanGrossPrasad}), in the case that $T$ is viewed as the unitary group of a one-dimensional Hermitian space over $K$.

\item It is the unique faithful irreducible representation of ${}^L T$, and even (among the irreducible representations) the unique one for which the counting problem is well-defined: for any other irreducible representation, there will be infinitely many characters of bounded analytic conductor.
\end{itemize}

\subsection{The Analytic Conductor}\label{theac}
The results of the preceding sections may be summarized by the following:
\begin{prop}\label{sharpmap}
For $\psi \in \cA(T)$, with the choice of the standard representation $\rho$ of the $L$-group of $T$, the $L$-function attached to $\psi$ is given by
$$
L(s, \psi, \rho) = L(s, \psi^\sharp)
$$ 
where $\psi^\sharp: \A_K^\times / K^\times \to \C^\times$ is the Hecke character $\psi^\sharp(x) = \psi(x/x^\sigma)$ and the right-hand side is the usual Hecke $L$-function.
\end{prop}

Given $\psi \in \cA(T)$, let $\fa$ be the largest ideal of $\cO_K$ such that $\psi^\sharp$ is a Hecke character modulo $\fa$. Then $\psi^\sharp$ has completed $L$-function \est{\Lambda(s,\psi^\sharp) = (|d_K|\N(\fa))^{s/2} \Gamma_{\C}(s + \kappa)L(s, \psi^\sharp),} where $d_K$ is the discriminant of the extension $K/\Q$ and $\kappa$ is a non-negative integer.  This completed $L$-function admits a functional equation \est{ \Lambda(s,\psi^\sharp) = W(\psi^\sharp) \Lambda(1-s,\overline{\psi^\sharp})} with $|W(\psi^\sharp)|=1$.  Note that $\overline{\psi^\sharp(x)} = {\psi^{\sharp}(x)}^{-1} = \psi^\sharp(x^\sigma)$ so that $L(s,\psi, \rho) =L(s,\overline{\psi},\rho)$ and, so in fact we have $W(\psi^\sharp)=\pm 1$ for any $\psi \in \cA(T)$, see for example \cite[\textsection 3.8]{IK} or \cite[Lecture 2]{RohrlichRootNumbers}.  
\begin{defin}[Analytic conductor with respect to $V_1$]
The positive real-valued function of $s\in \C$ given by \est{C(s,\psi) \eqdef |d_K|\N(\fa) (|s + \kappa| + 1)^{2}} is called the \emph{analytic conductor} of $\psi \in \cA(T)$.
\end{defin}
A few remarks:
\begin{enumerate}[\upshape (i)]
\item the normalization at the infinite place in the above definition differs from others in the literature in that $(|s + \kappa_j| + 1)$ is replaced by:
\begin{itemize}
\item $|s + \kappa_j| + 2$ (in \cite{MichelVenkateshGL2})
\item $|s + \kappa_j| + 3$ (in \cite{IK})
\item $|i \text{Im}(s) + \kappa_j| + 1$ (in \cite{IS2000GAFA}).
\end{itemize}
\item For counting purposes, we take $C(\psi)= C(0,\psi)$ and also call this the analytic conductor of $\psi$.
\item It is interesting to see what definition of analytic conductor one would obtain upon choosing a non-standard representation $\rho$ of the $L$-group.  If one takes $\rho = \C$ then $C(\psi)= 1$ for all $\psi \in \cA(T)$ and if one chooses $\rho = \C(-)$ one obtains $C(\psi) = |d_K|,$ or $4|d_K|$ or all $\psi \in \cA(T)$ depending on whether $\legen{d_K}{-1}$ is $+1$ or $-1$, respectively.  For either of these two choices the ``analytic conductor'' does not satisfy the finiteness property of height mentioned in the introduction.  For $V_n$ one gets the same definition for ``analytic conductor'' where the $\N(\fa)$ and $\kappa_j$ are now taken from the completed $L$-function for the Hecke character $(\psi^\sharp)^n$.  
\end{enumerate}

\section{The set-up for counting}\label{setup}

In this section we formulate the counting problem in classical terms which are more amenable to computation.  As in the previous section, write $S = \Res_{K/\Q} \G_m$ and $T = \ker \text{norm}: S \to \G_m$. Recall the map $r: S \to T$ given on $A$-points, for $A$ a $\Q$-algebra, by $x \mapsto x/x^\sigma$. By Hilbert's theorem 90, this map is surjective on $\Q$-points, so one has an exact sequence
$$
1 \to \G_m(\Q) \to S(\Q) \to T(\Q) \to 1.
$$
For any place $v$ of $\Q$, one similarly has
$$
1 \to \Q_v^\times \to (K \otimes \Q_v)^\times \to T(\Q_v) \to 1,
$$
where at the archimedean place and non-split places one again appeals to Hilbert's theorem 90, and at split places, the result is clear from the description of the norm map in \textsection\ref{torusfacts}.

The torus $S$ admits the integral model $\Res_{\cO_K/\Z} \G_m$, and extending the norm map in the obvious way to this model gives a model for $T$ as well; abusing notation, we will continue to write $S$ and $T$ for these models (note that they are not tori over $\text{Spec } \Z$). For $p$ a split prime, the argument of \textsection \ref{torusfacts}, with $\Q_p$ replaced by $\Z_p$, identifies $T(\Z_p)$ with the antidiagonal elements in $\Z_p^\times \times \Z_p^\times$ (if $p = 2$, the argument needs a modification, which we omit as we do not use it).

For all but finitely many primes, one has an exact sequence
$$
1 \to \G_m(\Z_p) \to S(\Z_p) \to T(\Z_p) \to 1.
$$
This is again obvious at split places (at least other than $p = 2$). At inert places, given $\beta \in T(\Z_p)$ there exists $\beta \in S(\Q_p)$ with $\beta^\sigma/\beta = \alpha$, and, as $p$ is a uniformizer of $K \otimes \Q_p$, one may multiply $\beta$ by a power of $p$ to make it a unit. (The sequence is not exact at ramified places).

It follows that there is an exact sequence
$$
1 \to \G_m(\A_\Q) \to S(\A_\Q) \to T(\A_\Q) \to 1.
$$

By Pontryagin duality, the sequences of character groups
$$
1 \to T(\Q)^\vee \to {K^\times}^\vee \to {\Q^\times}^\vee \to 1
$$
and
$$
1 \to T(\A_\Q)^\vee \to{\A_K^\times}^\vee \to {\A_\Q^\times}^\vee \to 1
$$
are also exact. As the restriction ${T(\A_\Q)}^\vee \to {T(\Q)}^\vee$ is surjective, one has an exact sequence
$$
1 \to \mathcal{A}(T) \to \mathcal{A}(S) \to \mathcal{A}(\G_m) \to 1.
$$
where for a commutative reductive group $G$, $\mathcal{A}(G)$ denotes the group of characters of $G(\A_\Q)$ trivial on $G(\Q)$. Caution: note that $\psi^\sharp|_{T(\A_f)} = \psi^2$.  

Now we write the characters $\psi^\sharp$ in terms of more computable classical data. The following notation system will be used frequently: if $A$ and $B$ are abelian subgroups of some ambient group, then $A^\vee \times_c B^\vee$ denotes the set of pairs $(\chi_A, \chi_B)$ of characters of $A$ and $B$ agreeing on $A \cap B$. (The subscript $c$ is for ``compatible''.)

Recall the exact sequence
$$
0 \to K^\times{\widehat{\cO}_K}^\times\C^\times \to \A_K^\times \to \text{Cl}(K) \to 0.
$$
Dualizing, one has
$$
0 \to \text{Cl}(K)^\vee \to  {\A_K^\times}^\vee \stackrel{\res}{\to} (K^\times{\widehat{\cO}_K}^\times\C^\times)^\vee \to 0
$$
The subgroup $\cA(S)$ of ${\A_K^\times}^\vee$ maps under $\res$ to the group of functionals trivial on $K^\times$, which is naturally isomorphic to
$$
{\C^\times}^\vee \times_c \left({\widehat{\cO}_K}^{\times}\right)^\vee,
$$
i.e.~to pairs $(\chi_\infty, \chi_f)$ of characters on $\C^\times$ and ${\widehat{\cO}_K}^\times$ whose restrictions to $\cO_K^\times$ agree. For $\chi$ a Hecke character of $K$, write $\res(\chi) = (\chi_\infty, \chi_f)$ under the above decomposition. 

The analytic conductor of $\chi$ (viewed as a form on $S$ with the standard representation $V_{1,0}$ of ${}^L S$ tacitly chosen, or equivalently viewed as a form on $\GL_1$ over $K$) is completely determined by the pair $(\chi_\infty, \chi_f)$, as follows from the calculation of the functional equation for the $L$-function of a Hecke character, as e.g.~in Tate's thesis. The local $L$-function of a $\chi_\infty \in {\C^\times}^\vee$ is of the form \est{ L(s,\chi_\infty) = \Gamma_\C(s+\kappa)} for some $\kappa \in \C$.  Then $c(\chi_\infty) = (|\kappa|+1)^2$ (following the conventions of \textsection \ref{theac} at the archimedean place).  The conductor ideal or ``finite conductor'' of $\chi$ is the largest ideal $I$ (smallest by divisibility) such that $\chi_f$ factors through $(\cO_K/I)^\times$. We write $I =c_f(\chi)$ or sometimes $I = c_f(\chi_f)$ for this conductor.  

Define the analytic conductor $C((\chi_\infty,\chi_f))$ of a pair such that for $\chi \in \cA(S)$, one has $C(\res(\chi)) = C(\chi,0)$, where $C(\chi, s)$ is the analytic conductor for $\chi$ as a form on $S$.  Explicitly, this amounts to \est{C((\chi_\infty,\chi_f)) = |d_K| \N\left(c_f(\chi_f)\right) c(\chi_\infty).} 

\begin{prop}\label{TheThingWeCount}
Automorphic forms $\phi \in \cA(T)$ are in analytic conductor-preserving $h_K$-to-one correspondence with pairs $(\chi_\infty,\chi_f) \in (\C^\times)^\vee \times_c (\widehat{\mathcal{O}}_K^\times)^\vee$ such that $\chi_\infty|_{\R^\times} = 1$ and $\chi_f|_{\widehat{\Z}^\times} = 1$. 

In particular, for a fixed positive integer $C$, the number of automorphic forms on $T$ of conductor $C$ is given by
\begin{multline*}  h_K \cdot \#\{(\chi_\infty,\chi_f) \in (\C^\times)^\vee \times_c (\widehat{\mathcal{O}}_K^\times)^\vee : \chi_\infty|_{\R^\times} = 1, \chi_f|_{\widehat{\Z}^\times} = 1, \\  \text{ and } c(\psi_\infty) \N(c_f(\chi_f))|d_K| = C \}.\end{multline*}
\end{prop}
\begin{proof}
The following diagram commutes ,and the vertical and horizontal sequences are short exact:
 \est{ \xymatrix{ & \operatorname{Cl}(K)^\vee \ar@{^{(}->}[d] &  \\ \mathcal{A}(T) \ar@{^{(}->}^\sharp[r] & \mathcal{A}(\G_m/K) \ar@{->>}[r]  \ar@{->>}[d]  & \mathcal{A}(\G_m/\mathbb{Q}) \bijar[d] \\ & (\C^\times)^\vee \times_c ( \widehat{\mathcal{O}}^\times)^\vee \ar@{->>}[r] & (\R^\times)^\vee \times_c (\widehat{\Z}^\times)^\vee }} The diagram characterizes the image of the map $\sharp$ in the space of Hecke characters in terms of the constituents $(\chi_\infty, \chi_f)$. For each pair $(\chi_\infty, \chi_f)$ satisfying the conditions of the proposition (i.e.~in the bottom middle of the diagram), there are exactly $h_K$ Hecke characters of $K$ lifting it. Each such lift has analytic conductor $c(\chi_\infty)\N(c_f(\chi_f))|d_K|$ and is of the form $\phi^\sharp$ for a unique $\phi \in \mathcal{A}(T)$, from which the proposition follows.
 \end{proof}

\section{The Generating Series}\label{GeneratingSeries}

It remains to count the set in Proposition \ref{TheThingWeCount}.  For $n \geq 1$, let $\Phi(n)$ be \begin{multline*} \Phi(n) \eqdef \#\{(\chi_\infty,\chi_f) \in (\C^\times)^\vee \times_c (\widehat{\mathcal{O}}^\times)^\vee : \chi_\infty|_{\R^\times} = 1, \chi_f|_{\widehat{\Z}^\times} = 1, \\ \text{ and } c(\psi_\infty) \mathbf{N}(c_f(\chi_f)) = n \}\end{multline*} and \est{\Phi_1(n) \eqdef \#\{ \chi \in   (\widehat{\mathcal{O}}^\times)^\vee : (1,\chi) \in (\C^\times)^\vee \times_c (\widehat{\mathcal{O}}^\times)^\vee, \chi |_{\widehat{\Z}^\times} = 1, \mathbf{N}(c_f(\chi)) = n \}.}  By Proposition \ref{TheThingWeCount} one has (cf.~Theorem \ref{MT}): \es{\label{mt1form1}  \sum_{\substack{\psi \in \cA(T) \\ C(\psi) \leq X}} 1 = h_K \sum_{|d_k| n \leq X} \Phi(n) } and \es{\label{mt1form2} \sum_{\substack{\psi \in \cA(T) \\ \psi(K_\infty^\times ) \\ C(\psi) \leq X}} 1= h_K \sum_{|d_k| n \leq X} \Phi_1(n) .}
Let \est{\phi^*(\mathfrak{a}) = \# \{ \chi_f \in (\widehat{\mathcal{O}}^\times)^\vee : \chi_f|_{\widehat{\Z}^\times} = 1, \text{ and } c_f(\chi_f)=\mathfrak{a}\}} and \est{  \phi_\infty(n) = \#\{ \chi_\infty \in (\C^\times)^\vee : \chi_\infty|_{\R^\times} = 1,\text{ and } c(\chi_\infty)=n\}.}

In the case that $K\neq \Q(i),\Q(\zeta_3)$ the compatibility condition $\times_c$ is automatically satisfied, i.e.~ $\times_c = \times$.  For $K\neq \Q(i),\Q(\zeta_3)$, it follows that \es{\label{P1} \Phi(n)= \sum_{ab=n} \left(\sum_{\N(\mathfrak{a}) = a} \phi^*(\mathfrak{a}) \right) \phi_\infty(b)} and \es{\label{P11} \Phi_1(n) =  \sum_{\N(\mathfrak{a}) = n} \phi^*(\mathfrak{a}).}  
In the case that $K=\Q(i)$, set for $u=1,-1$ \est{\phi^*_u(\mathfrak{a}) = \#\{ \chi \in (\widehat{\mathcal{O}}^\times)^\vee : \chi|_{\widehat{\Z}^\times} = 1, c_f(\chi)=\mathfrak{a}, \text{ and } \chi(i) = u\}} and similarly, in the case $K=\Q(\zeta_3)$, set $u=1,\zeta_3,\overline{\zeta_3}$ we set \est{\phi^*_u(\mathfrak{a}) = \#\{ \chi \in (\widehat{\mathcal{O}}^\times)^\vee : \chi|_{\widehat{\Z}^\times} = 1, c_f(\chi)=\mathfrak{a}, \text{ and } \chi(\zeta_6) = u\}.}  There is a conflict of notation if $u=1$, but in this case it will be clear from context which of the above two functions is intended by $\phi^*_u(\mathfrak{a})$.  
In the case $K=\Q(i)$, \es{\label{P2} \Phi(n) = \sum_{ab=n} \left(\sum_{\N(\mathfrak{a}) = a} \phi^*_{(-1)^{b-1}}(\mathfrak{a}) \right) \phi_\infty(b)} and \es{\label{P12}  \Phi_1(n) =  \sum_{\N(\mathfrak{a}) = n} \phi^*_1(\mathfrak{a}).}  
In the case $K=\Q(\zeta_3)$, \es{\label{P3} \Phi(n) = \sum_{ab=n} \left(\sum_{\N(\mathfrak{a}) = a} \phi^*_{\zeta_3^{\sqrt{b}-1}}(\mathfrak{a}) \right) \phi_\infty(b)} and \es{\label{P13}  \Phi_1(n) =  \sum_{\N(\mathfrak{a}) = n} \phi^*_1(\mathfrak{a}).}  

The count at archimedean places follows from the following lemma.
\begin{lem}\label{archphi} We have
\est{ \phi_\infty(n) =  \begin{cases} 2 & \text{ if } n\geq 4, n= \square \\ 1 & \text{ if } n=1 \\ 0 & \text{ else.}\end{cases}.}
\end{lem}
\begin{proof}
The characters of $T(\R)=S^1$ are of the form $\chi_k(z)= z^k$ for $k\in \Z$, and thus $\chi_k^\sharp = z^{2k}|z|^{-2k}$ with $k \in \Z$.  So, for $\psi_\infty \in T(\R)$, $\psi_\infty(z) = z^{2k}|z|^{-2k}$ with $k \in \Z$; write $k=k_\psi$ for this integer.  Recall from  \textsection \ref{setup} that $c(\psi_\infty) = (1+|k_\psi|)^2.$  
\end{proof}
We now turn to the non-archimedean places.  Let  \est{ \phi(\mathfrak{a}) = \#\{ \chi_f \in (\widehat{\mathcal{O}}^\times)^\vee : \chi_f|_{\widehat{\Z}^\times} = 1, \text{ and } \chi_f \text{ factors through } (\mathcal{O}/\mathfrak{a})^\times\}. }  By definition, $\chi_f$ factors through $(\mathcal{O}/\mathfrak{a})^\times$ if and only if $c_f(\chi_f) \mid \mathfrak{a}$, so \est{\sum_{\mathfrak{d} \mid \mathfrak{a}} \phi^*(\mathfrak{d}) = \phi(\mathfrak{a}).}  Let $\mu(\mathfrak{d})$ be the M\"obius function on ideals of $\mathcal{O}$. By M\"obius inversion, \es{\label{oldLem2} \phi^*(\mathfrak{a}) = \sum_{\mathfrak{d} \mid \mathfrak{a}} \mu(\mathfrak{d}) \phi(\mathfrak{a}\mathfrak{d}^{-1}).}

\begin{lem}\label{nonprimitivecharacters}
Let $\mathfrak{p}, \overline{\mathfrak{p}}$ be a pair of primes of $\mathcal{O}$ lying over a split prime $p$. Then
 \est{\phi(\mathfrak{p}^{n} \overline{\mathfrak{p}}^m ) = \begin{cases} 1 & \text{ if } \min(n,m)=0 \\ p^{\min(n,m)-1}(p-1) & \text{ else.}\end{cases}} 
If $\mathfrak{p}$ is lying over an inert prime $p$, then
\est{\phi(\mathfrak{p}^n) = \begin{cases} 1 & \text{ if } n=0 \\p^{n-1}(p+1) & \text{else,} \end{cases}} and if $\mathfrak{p}^2$ is lying over a ramified prime $p$ ,then \est{ \phi(\mathfrak{p}^n) = p^{\lfloor \frac{n}{2} \rfloor }.}
\end{lem}

\begin{proof} For an ideal $\mathfrak{a}$, the quantity $\phi(\mathfrak{a})$ is equal to the number of characters of $\left(\mathcal{O}_K/\mathfrak{a}\right)^\times$ which are trivial on $\left(\Z/\mathfrak{a}\cap \Z\right)^\times$, which is just equal to the index of the finite subgroup $\left(\Z/\mathfrak{a}\cap \Z \right)^\times$ in $\left(\mathcal{O}_K/\mathfrak{a}\right)^\times$. The orders of these groups are readily computable in the case where $\mathfrak{a}$ is one of the ideals in the lem.
\end{proof}

The next argument moves from $\phi$ to $\phi^*$ with Mobius inversion, i.e.~\eqref{oldLem2}.

\begin{lem}\label{philem}
We have if $\mathfrak{p}, \overline{\mathfrak{p}}$ be a pair of primes of $\mathcal{O}$ lying over a split prime $p$ then 
\est{\phi^*(\mathfrak{p}^{n} \overline{\mathfrak{p}}^m ) = \begin{cases}   p -2 & \text{ if } n=m=1 \\ p^{n-2}(p-1)^2 & \text{ if } n=m\geq 2 \\ 0 & \text{ if } n \neq m,\end{cases}} and 
\est{ \sum_{\N(\mathfrak{a})=p^n}\phi^*(\mathfrak{a}) = \begin{cases} p- 2& \text{ if } n=2 \\ p^{\frac{n}{2}-2}(p-1)^2 & \text{ if } n \geq 4 \text{ even} \\  0 & \text{ otherwise}.\end{cases}}   
If $\mathfrak{p}$ lies over an inert prime $p$ we have  \est{\phi^*(\mathfrak{p}^n) = \begin{cases} p & \text{ if } n=1 \\ p^n-p^{n-2} & \text{ if } n\geq 2, \end{cases}} and \est{ \sum_{\N(\mathfrak{a})=p^n}\phi^*(\mathfrak{a}) = \begin{cases} p & \text{ if } n=2 \\ p^{\frac{n}{2}} -p^{\frac{n}{2}-2} & \text{ if } n \geq 4 \text{ even} \\  0 & \text{ otherwise}.\end{cases}} If $\mathfrak{p}$ lies over a ramified prime $p$ then \est{ \phi^*(\mathfrak{p}^n) =\begin{cases} p^{\frac{n}{2}}-p^{\frac{n}{2}-1} & \text{ if } n \text{ even} \\ 0 & \text{ if } n \text{ odd},\end{cases}} and \est{ \sum_{\N(\mathfrak{a})=p^n}\phi^*(\mathfrak{a}) = \begin{cases} p^\frac{n}{2}-p^{\frac{n}{2}-1} & \text{ if } n \text{ even} \\ 0 & \text{ if } n \text{ odd}.\end{cases}}
\end{lem}
\begin{proof}
The statements in the cases of $p$ inert and ramified are a simple consequence of \eqref{oldLem2} and Lemma \ref{nonprimitivecharacters} along with the fact that $N(\mathfrak{p})=p^2$ if $\mathfrak{p}$ lies over an inert prime $p$, and $N(\mathfrak{p})=p$ if $\mathfrak{p}$ lies over a ramified prime $p$.

Now consider the case that $\mathfrak{p},\overline{\mathfrak{p}}$ are a pair of primes lying over a split prime $p$.  If $n\neq m$, we may assume without loss of generality that $n>m$ and in this case by \eqref{oldLem2}, \est{ \phi^*(\mathfrak{p}^n\overline{\mathfrak{p}}^m) =  \phi(\mathfrak{p}^n\overline{\mathfrak{p}}^m) -\phi(\mathfrak{p}^{n-1}\overline{\mathfrak{p}}^m) - \left( \phi(\mathfrak{p}^n\overline{\mathfrak{p}}^{m-1}) -\phi(\mathfrak{p}^{n-1}\overline{\mathfrak{p}}^{m-1})\right) =0} by Lemma \ref{nonprimitivecharacters} since in any case $n-1\geq m$.  

If $n=m$, then $\phi^*(\mathfrak{p}^n\overline{\mathfrak{p}}^n)$ can be calculated directly from \eqref{oldLem2} and Lemma \ref{nonprimitivecharacters}.  There is at most one non-vanishing summand in $\sum_{N(\mathfrak{a})=p^n}\phi^*(\mathfrak{a})$ if $p$ is split because $\phi^*(\mathfrak{p}^n\overline{\mathfrak{p}}^m)$ vanishes unless $n=m$.  If $n=m$ is even, then there is one non-vanishing term, and if $n=m$ is odd, there are none.
\end{proof}

This proves Theorem \ref{MT2} in the case that $K \neq \Q(i),\Q(\zeta_3)$.  We now give the analogous result for the functions $\phi_u^*(\mathfrak{a})$ for $\Q(i)$ and $\Q(\zeta_3)$ to finish the proof of Theorem \ref{MT2}.

\begin{lem}\label{phiu}
If $K=\Q(i)$ and $u=\pm1$, then \est{ \phi^*_u(\fa) = \frac{1}{2} \phi^*(\fa) +u\frac{1}{2} \sum_{\substack{\mathfrak{d} \mid \fa \\ (n) \nmid \mathfrak{d} \\ \text{ for all } n}}  \mu(\fa/\mathfrak{d}).} If $K=\Q(\zeta_3)$ and $u=1,\zeta_3,\overline{\zeta_3}$, then \est{ \phi^*_u(\fa) = \frac{1}{3} \phi^*(\fa) +(u+\overline{u})\frac{1}{3} \sum_{\substack{\mathfrak{d} \mid \fa \\ (n) \nmid \mathfrak{d} \\ \text{ for all } n}}  \mu(\fa/\mathfrak{d}).}
\end{lem}
\begin{proof}
Let $\delta_A$ be $0$ or $1$ according as the statement $A$ is false or true.  The condition that $\chi(i)$ or $\chi(\zeta_6)$ be a specified root of unity can be expressed via the following formulas:  for $K=\Q(i)$ and $u=\pm 1$, \est{ \delta_{\chi(i)=u} = \frac{1}{2} (1+u \chi(i)),} and for $K=\Q(\zeta_3)$ and $u=1,\zeta_3,\overline{\zeta_3}$, one has \est{ \delta_{\chi(\zeta_6) = u} = \frac{1}{3} (1+ \overline{u} \chi(\zeta_6) + u \chi(\zeta_3)).}  Thus for $\Q(i)$, \est{ \phi^*_u(\mathfrak{a}) = & \frac{1}{2}\sum_{\substack{c_f(\chi) = \mathfrak{a} \\ \chi |_{\widehat{\Z}^\times} =1}} 1 + \frac{1}{2}\sum_{\substack{c_f(\chi) = \mathfrak{a} \\ \chi |_{\widehat{\Z}^\times} =1}}\chi(i) \\ = & \frac{1}{2} \phi^*(\mathfrak{a}) + \frac{1}{2}\sum_{\substack{c_f(\chi) = \mathfrak{a} \\ \chi |_{\widehat{\Z}^\times} =1}}\chi(i),} and similarly for $\Q(\zeta_3)$, \est{ \phi^*_u(\mathfrak{a}) = \frac{1}{3} \phi^*(\mathfrak{a}) + \frac{1}{3}\sum_{\substack{c_f(\chi) = \mathfrak{a} \\ \chi |_{\widehat{\Z}^\times} =1}}\chi(\zeta_6)+ \frac{1}{3}\sum_{\substack{c_f(\chi) = \mathfrak{a} \\ \chi |_{\widehat{\Z}^\times} =1}}\chi(\zeta_3).}  These last sums can be calculated using orthogonality of characters and M\"obius inversion.    

Indeed, for $x\in \widehat{\mathcal{O}}^\times$, \begin{multline*} \sum_{\substack{c_f(\chi) \mid \mathfrak{a} \\ \chi |_{\widehat{\Z}^\times}=1}} \chi(x) \\ = \begin{cases} \phi(\mathfrak{a}) & \text{ if } x \in \langle \widehat{\Z}^\times, 1 \pmod*{\mathfrak{a}}\rangle = \text{Image}((\Z/\mathfrak{a} \cap \Z)^\times \hookrightarrow (\cO / \mathfrak{a})^\times) \\ 0 & \text{ else}\end{cases}\end{multline*} by orthogonality of characters.  Writing $I_\mathfrak{a}$ for the image of $(\Z/\mathfrak{a} \cap \Z)^\times$ in $(\cO/\mathfrak{a})^\times$, it suffices to determine for which ideals of $\mathbb{Z}[i]$ one has  $i \in I_\fa$ or similarly for which ideals of $\mathbb{Z}[\zeta_3]$ one has $\zeta_6,\zeta_3 \in I_\fa$. 
 
We claim that any of $i,\zeta_6,\zeta_3 \in I_\fa$ if and only if $\phi(\mathfrak{a})=1$ if and only if $(p) \nmid \fa$ for all rational primes $p$.  The second equivalence in this claim can be seen immediately from the computation in Lemma \ref{nonprimitivecharacters}.  

For the first equivalence, in the 	``if'' direction, note that if $\phi(\mathfrak{a})=1$ then $i,\zeta_6,\zeta_3 \in (\cO/\fa)^\times= (\Z/\fa \cap \Z)^\times$ by the proof of Lemma \ref{nonprimitivecharacters}. Conversely, if $\mathfrak{b} \mid \fa$ and e.g.~$\zeta_6 \in I_\fa$, then automatically $\zeta_6 \in I_\fb.$  Thus it suffices to show that for each rational prime $p$ one has $i,\zeta_6,\zeta_3 \not \in I_{(p)},$ but this is clear.

Therefore, by M\"obius inversion again, for $u=i,\zeta_6,\zeta_3$ one has \est{\sum_{\substack{c_f(\chi) = \mathfrak{a} \\ \chi |_{\widehat{\Z}^\times} =1}}\chi(u) = \sum_{\substack{\mathfrak{d} \mid \fa \\ (p) \nmid \mathfrak{d} \\ \text{ for all } p}} \phi(\mathfrak{d}) \mu(\fa/\mathfrak{d}).}  By the claim, $(p) \nmid \mathfrak{d}$ for all $p$ if and only if $\phi(\mathfrak{d})=1$, so in fact \est{\sum_{\substack{c_f(\chi) = \mathfrak{a} \\ \chi |_{\widehat{\Z}^\times} =1}}\chi(u) = \sum_{\substack{\mathfrak{d} \mid \fa \\ (p) \nmid \mathfrak{d} \\ \text{ for all } p}} \mu(\fa/\mathfrak{d}).} 
\end{proof}

Now for $\real(s)>1$, let \est{ L(s,\Phi) \eqdef \sum_{n\geq 1} \frac{\Phi(n)}{n^s}, \,\,\,\,\,\,\, L(s,\Phi_1) \eqdef \sum_{n\geq 1} \frac{\Phi_1(n)}{n^s},}\est{ L(s,\phi^*) \eqdef \sum_{\mathfrak{a}} \frac{\phi^*(\mathfrak{a})}{\N(\mathfrak{a})^s}, \,\,\, \text{ and }\,\,\, L(s,\phi_u^*) \eqdef \sum_{\mathfrak{a}} \frac{\phi^*_u(\mathfrak{a})}{\N(\mathfrak{a})^s}.} 
The relations among these generating series are summarized in the following lemma:
\begin{lem}\label{genseriesrlns}
 When $K \neq \Q(i),\Q(\zeta_3)$, \est{ L(s,\Phi)  = (2 \zeta(2s)-1)L(s,\phi^*)} and \est{ L(s,\Phi_1) = L(s,\phi^*).}  When $K=\Q(i)$, \est{ L(s,\Phi) = 2 \left(1-\frac{1}{2^{2s}}\right) \zeta(2s)L(s,\phi_1^*) + \frac{2}{2^{2s}} \zeta(2s) L(s,\phi^*_{-1}) - L(s,\phi_1^*)} and \est{L(s,\Phi_1) = L(s,\phi_1^*).}  When $K=\Q(\zeta_3)$, \begin{multline*} L(s,\Phi) =   \left(  \left(1-\frac{1}{3^{2s}}\right)\zeta(2s) + L(2s,\chi_3)\right) L(s,\phi^*_1) \\  +  \left(  \left(1-\frac{1}{3^{2s}}\right)\zeta(2s) - L(2s,\chi_3)\right) L(s,\phi^*_{\zeta_3}) \\  + 2 \frac{1}{3^{2s}}\zeta(2s)L(s,\phi^*_{\overline{\zeta_3}}) - L(s,\phi^*_1),\end{multline*} where $L(s,\chi_3)$ is the classical Dirichlet $L$-function of the unique non-trivial Dirichlet character mod $3$, and \est{L(s,\Phi_1) = L(s,\phi_1^*).} \end{lem}
\begin{proof}
The first, second, fourth and sixth statements of the lemma follow directly from equations \eqref{P1}, \eqref{P11}, \eqref{P12}, \eqref{P13} and Lemma \ref{archphi}.  

For the third statement, note from \eqref{P2} and Lemma \ref{archphi} that if $K=\Q(i)$ then \begin{multline*} L(s,\Phi) =  2 \left( \sum_{n\equiv 1 \pmod*{2}}\frac{1}{n^{2s}} \right) L(s,\phi^*_1) \\  + 2 \left( \sum_{n\equiv 0 \pmod*{2}}\frac{1}{n^{2s}} \right) L(s,\phi^*_{-1}) - L(s,\phi^*_1)\end{multline*} and if $K=\Q(\zeta_3)$ then \begin{multline*} L(s,\Phi) =  2 \left( \sum_{n\equiv 1 \pmod*{3}}\frac{1}{n^{2s}} \right) L(s,\phi^*_1) \\  + 2 \left( \sum_{n\equiv 2 \pmod*{3}}\frac{1}{n^{2s}} \right) L(s,\phi^*_{\zeta_3}) \\  + 2 \left( \sum_{n\equiv 0 \pmod*{3}}\frac{1}{n^{2s}} \right) L(s,\phi^*_{\overline{\zeta_3}}) - L(s,\phi^*_1).\end{multline*} The lemma then follows from standard simplifications of the above Dirichlet series.  
\end{proof}
In the special cases $K=\Q(i)$ and $K=\Q(\zeta_3)$, we additionally relate the functions $L(s,\phi^*_u)$ back to $L(s,\phi^*)$.
\begin{lem}
If $K=\Q(i)$, \est{L(s,\phi^*_{\pm1}) = \frac{1}{2} L(s,\phi^*) \pm \frac{1}{2} \zeta(2s)^{-1}.} When $K=\Q(\zeta_3)$, \est{ L(s,\phi_1^*) = \frac{1}{3} L(s,\phi^*) + \frac{2}{3} \zeta(2s)^{-1}} and \est{L(s,\phi_{\zeta_3}^*) = L(s,\phi_{\overline{\zeta_3}}^*) = \frac{1}{3} L(s,\phi^*) - \frac{1}{3} \zeta(2s)^{-1}.}
\end{lem}
\begin{proof} Let \est{W_K(s) = \sum_{\substack{\fa \\ (n) \nmid \fa \\ \text{for all } n}} \frac{1}{N(\fa)^s}.}
Directly from Lemma \ref{phiu}, for $K=\Q(i)$, \est{L(s,\phi^*_{\pm1}) = \frac{1}{2} L(s,\phi^*) \pm \frac{1}{2} \frac{1}{\zeta_K(s)} W_K(s),} and for $K=\Q(\zeta_3)$, \est{ L(s,\phi_1^*) = \frac{1}{3} L(s,\phi^*) + \frac{2}{3} \frac{1}{\zeta_K(s)} W_K(s)} and \est{L(s,\phi_{\zeta_3}^*) = L(s,\phi_{\overline{\zeta_3}}^*) = \frac{1}{3} L(s,\phi^*) - \frac{1}{3}\frac{1}{\zeta_K(s)}  W_K(s).}
For all ideals $\fa$ there exists a unique pair $n\in \mathbb{N}$ and $\mathfrak{s}$ having the property that $(m) \nmid \mathfrak{s}$ for any $m\in \mathbb{N}$ such that $\fa=(n)\mathfrak{s}$.  Then \est{ \zeta_K(s) = \sum_{\fa} \frac{1}{N(\fa)^s} = \left( \sum_{n\geq 1} \frac{1}{n^{2s}} \right) \left( \sum_{\substack{\mathfrak{s} \\ (m) \nmid \mathfrak{s} \\ \text{for all } m}} \frac{1}{N(\mathfrak{s})^s} \right) = \zeta(2s) W_K(s).}  Thus $W_K(s) = \zeta_K(s)\zeta(2s)^{-1}.$
\end{proof}

\begin{prop}\label{phi*prop}
One has \es{\label{phi*formula}L(s,\phi^*) = \frac{\zeta(2s-1)}{\zeta_K(2s)}.}
\end{prop}
\begin{proof}
Consider the case of $p$ a split prime.  Expanding in geometric series, \begin{multline*}\left( 1-\frac{1}{p^{2s}}\right)^2\left(1-\frac{1}{p^{2s-1}}\right)^{-1} \\ = 1+ \frac{p-2}{p^{2s}} + \frac{p^2-2p+1}{p^{4s}} + \frac{p^3-2p^2+p}{p^{6s}} + \cdots \end{multline*} which by Lemma \ref{philem} shows that the Euler factors at split primes on both sides of \eqref{phi*formula} match.  
If $p$ is an inert prime, then \est{ \left( 1-\frac{1}{p^{4s}}\right)\left(1-\frac{1}{p^{2s-1}}\right)^{-1} = 1+\frac{p}{p^{2s}} + \frac{p^2-1}{p^{4s}} +\frac{p^3-p}{p^{6s}}+\cdots } and if $p$ is a ramified prime, then \est{ \left( 1-\frac{1}{p^{2s}}\right)\left(1-\frac{1}{p^{2s-1}}\right)^{-1} = 1+\frac{p-1}{p^{2s}} + \frac{p^2-p}{p^{4s}} + \cdots .}  These computations show that the Euler factors at inert and ramified primes match on both sides of \eqref{phi*formula}.
\end{proof}
 
Note in particular that for each fixed $\eps>0$ one has, uniformly in $K$ and $s$ in the region $\real(s)\geq 1/2+\eps$, that \es{\label{DedekindBound}\frac{1}{\zeta_K(2s)} \leq \zeta(1+2\eps)^2\leq (1 + (2\eps)^{-1})^2.}

We next use Perron's formula to move from the Dirichlet series just computed to the counts in Theorem \ref{MT}.  Perron's formula (see e.g.~\cite[Lemma page 105]{Davenport}) states that for any $\kappa > 0$, $T>0$, and $y>0$, \est{\frac{1}{2 \pi i } \int_{\kappa -iT}^{\kappa + iT} y^s\,\frac{ds}{s} = \begin{cases} 1 & \text{ if } y>1 \\ 0 & \text{ if }y<1\end{cases} + O \left( \frac{y^\kappa}{\max(1, T |\log y|)}\right).}  By Lemma \ref{philem}, $\Phi(n) \leq 2\delta_{n=\square}\sigma_1(\sqrt{n}) \ll \delta_{n= \square} \sqrt{n} \log \log n$ if $n\geq 3$.  To prove the first formula in Theorem \ref{MT} in the cases $K\neq \Q(i),\Q(\zeta_3)$ it suffices by \eqref{mt1form1}, Lemma \ref{genseriesrlns}, Proposition \ref{phi*prop} and Perron's formula to estimate for $Y\geq 3$ \begin{multline}\label{Perronstart} \sum_{n \leq Y} \Phi(n) = \frac{1}{2\pi i } \int_{\kappa-iT}^{\kappa+iT} \frac{(2\zeta(2s)-1)}{\zeta_K(2s)} \zeta(2s-1)Y^s \frac{ds}{s} \\ + O \left(Y \sum_{n=\square} \frac{\sigma_1(\sqrt{n})}{n^\kappa \max(1, T |\log Y/n|)} \right),\end{multline} where we have chosen $\kappa = 1+1/\log Y$, $T=Y^\alpha$, and $0<\alpha\leq 1$ to be determined later.    

We first estimate the $O$ term above.  Let $U=\exp((2T)^{-1})$.  Make the change of variable $m^2=n$ and break the range of summation into three parts: $m\leq Y^{1/2}/U$, $Y^{1/2}/U< m \leq Y^{1/2}U$, and $m> Y^{1/2}U$.  In the first range, the $O$ term in \eqref{Perronstart} is bounded by \begin{multline*} \ll \frac{Y \log \log Y}{T} \sum_{m =1}^{Y^{1/2}/U} \frac{1}{m \log (X/m^2)}  \ll \frac{Y \log \log Y}{T} \int_{1/T}^{\log Y} \frac{du}{u} \\ \ll \frac{Y \log Y \log \log Y}{T}\end{multline*} by a change of variables.  In the second range, the $O$ term in  \eqref{Perronstart} is bounded by \est{ \ll Y \log \log Y \sum_{m=Y^{1/2}/U+1}^{Y^{1/2}U} \frac{1}{m}  \ll \frac{Y \log \log Y}{T}.} In the third range, the $O$ term in  \eqref{Perronstart} is bounded by \begin{multline*} \ll \frac{Y}{T} \sum_{m > Y^{1/2}U} \frac{\log \log m}{m^{1+2/\log Y} \log (m^2/Y)} \\ \ll \frac{Y}{T} \int_{1/T}^\infty \frac{\log \left(\frac{v+\log Y}{2}\right) \exp(-v/\log Y)}{v} \,dv\end{multline*} by a change of variables.  But $\frac{v+\log Y}{2} \leq \max(v,\log Y)$, so upon splitting the integral and making a change of variables, one sees that the sum in this third piece is also $\ll Y \log Y/T$.  

The next estimate is for the contour integral term in \eqref{Perronstart}.  Let $\eps>0$. Shift the contour to the other three sides of a box having corners at $\kappa+iT, 1/2+\eps+iT,1/2+\eps-iT,\kappa-iT$ and pick up the residue at $s=1$.  Then \begin{multline}\label{aftercontour}\sum_{m \leq Y} \Phi(n)  =  \frac{1}{2} \frac{2\zeta(2)-1}{\zeta_K(2)} Y   \\ + \left( \int_{\kappa + iT}^{1/2+\eps+ iT} + \int_{1/2+\eps+iT}^{1/2+\eps-iT} + \int_{1/2+\eps -iT}^{\kappa -iT}\right) \frac{2 \zeta(2s)-1}{\zeta_K(2s)} \zeta(2s-1) Y^s \frac{ds}{s} \\  + O\left(\frac{Y \log Y \log \log Y}{T}\right).\end{multline}  Using the bound \eqref{DedekindBound} and a similar bound for the Riemann zeta function, the first of the three integrals in \eqref{aftercontour} is bounded by \est{\ll (1+(2\eps)^{-1})^3 \int_{1/2+\eps}^{1+1/\log Y} \frac{|\zeta(2u-1+2iT)|}{|u+iT|}Y^{u}\,du,} which is, by a convexity bound on the zeta function, \est{ \ll (1+(2\eps)^{-1})^3 \int_{1/2+\eps}^{1+1/\log Y} \left(\frac{Y}{T}\right)^{u}\,du \ll \frac{ (1+(2\eps)^{-1})^3}{\log Y} \frac{Y}{T}.}  The third of three integrals in \eqref{aftercontour} is treated similarly.  Similarly, the second integral in \eqref{aftercontour} is \est{ \ll (1+(2\eps)^{-1})^3 Y^{1/2+\eps} \int_{-T}^T \frac{|\zeta(2 \eps + 2i t)| }{|\frac{1}{2} + \eps +it|}\,dt \ll (1+(2\eps)^{-1})^3 Y^{1/2+\eps} T^{1/2}.}  

Grouping these estimates together, \begin{multline*}\sum_{n \leq Y} \Phi(n) = \frac{1}{2} \frac{2\zeta(2)-1}{\zeta_K(2)} Y + O \left( (1+(2\eps)^{-1})^3 Y^{1/2+\eps} T^{1/2}\right. \\ \left. + \frac{ (1+(2\eps)^{-1})^3}{\log Y} \frac{Y}{T}+ \frac{Y \log Y \log \log Y}{T} \right).\end{multline*}  Now taking $T=Y^{1/3}$, $ \eps = 1/\log Y$ and inserting this to \eqref{mt1form1} gives the first statement of Theorem \ref{MT} in the case $K \neq \Q(i),\Q(\zeta_3)$.  When $K= \Q(i)$, \est{ L(s,\Phi) = \frac{1}{2} \frac{2 \zeta(2s)-1}{\zeta_K(2s)} \zeta(2s-1) + \frac{1}{2} \left( 2(1-\frac{2}{2^{2s}}) \zeta(2s)-1\right)\zeta(2s)^{-1}, } and when $K=\Q(\zeta_3)$, \begin{multline*} L(s,\Phi) \\ = \frac{1}{3} \frac{2\zeta(2s)-1}{\zeta_K(2s)} \zeta(2s-1) + \left(\frac{1}{3} \left(1-\frac{3}{3^{2s}}\right) \zeta(2s) + L(2s,\chi)-\frac{2}{3}\right) \zeta(2s)^{-1}.\end{multline*}  Running the same computation as above with these Dirichlet series completes the proof of the first claim in Theorem \ref{MT}.

For the second claim in Theorem \ref{MT} we work with $L(s,\Phi_1).$  In the $K\neq \Q(i),\Q(\zeta_3)$ case, \est{L(s,\Phi_1) = \frac{1}{\zeta_K(2s)} \zeta(2s-1).} In the $K=\Q(i)$ case,    \est{ L(s,\Phi_1) = \frac{1}{2} \frac{1}{\zeta_K(2s)}\zeta(2s-1)+ \frac{1}{2} \zeta(2s)^{-1},} and in the $K=\Q(\zeta_3)$ case, \est{ L(s,\Phi_1) = \frac{1}{3} \frac{1}{\zeta_K(2s)}\zeta(2s-1)+ \frac{2}{3} \zeta(2s)^{-1}.}
Running the same computation as above with these Dirichlet series completes the proof of the second claim of Theorem \ref{MT}.

\appendix
\section{Explicit calculation of the Langlands parameter}
This appendix contains a proof of Lemma \ref{cocycle}. The lemma holds in a bit more generality; thus, for the appendix only, the extension $K/F$ is an arbitrary quadratic extension of number fields or local fields. In the local case, $C_K$ denotes $K^\times$, and in the global case, it denotes the id\`{e}le class group $\mathbb{A}_K^\times/K^\times$; $C_K^{(1)}$ denotes the norm one elements of $C_K$.

Before proving the lemma, we summarize the construction of \cite{LanglandsAbelian}; in this summary, $T$ denotes an arbitrary torus over $F$, and $K/F$ can be any extension splitting $T$ (not just a quadratic one).  Write $L$ for the character lattice of $T$, and $\widehat{L}$ for the cocharacter lattice. Write $T(F)^!$ for $T(F)$ if $F$ is local, and for $T(\A_F)/T(F)$ is $F$ is global. Then one has 

\begin{align}
T(F)^! &\to \Hom_{\Gal(K/F)}(L, C_K) \label{standard1}\\
&= (C_K \otimes \widehat{L})^{\Gal(K/F)} \label{standard2}  \\
&= H_1(C_K, \widehat{L})^{\Gal(K/F)} \label{standard3} \\
&= H_1(W_{K/F}, \widehat{L}) \label{group homology, restriction map}
\end{align}
Maps (\ref{standard1}) - (\ref{standard3}) are straightforward, whereas identification (\ref{group homology, restriction map}) follows from a computation of Langlands showing that the trace map
$$
H_1(W_{K/F}, \widehat{L}) \to H_1(C_K, \widehat{L})^{\Gal(K/F)}
$$
is an isomorphism (this map is called restriction in \cite{LanglandsAbelian}). Writing $\widehat{T}$ for the dual torus to $T$, the cup product gives a perfect pairing
\begin{equation}
H_1(W_{K/F}, \widehat{L}) \otimes H^1(W_{K/F}, \widehat{T}) \to \C^\times, \label{cup product}
\end{equation}
Dualizing, there is therefore a map
\begin{equation}\label{langlandscorrespondence}
H^1(W_{K/F}, \widehat{T}) \to {T(F)^!}^{\vee}.
\end{equation} 
In the local case, map (\ref{standard1}), and hence map (\ref{langlandscorrespondence}), is a bijection and the desired correspondence; in the global case, (\ref{langlandscorrespondence}) is a finite-to-one surjection, but any choice of $\xi \in H^1(W_{K/F}, \widehat{T})$ mapping to a fixed $\psi \in {T(F)^!}^\vee$ gives rise to an equivalent Langlands parameter.

Returning to the case where $T$ is the non-split torus associated to $K/F$, choosing a lift $\widetilde{\sigma}$ of the generator $\sigma$ of $\Gal(K/F)$ to $W_{K/F}$, and given a character $\psi$ of $T(F)$, Lemma \ref{cocycle} asserts the following of the Langlands parameter $\alpha_\psi$:
\begin{itemize}
\item It is trivial on the kernel $[W_K, W_K]$ of the map $W_F \to W_{K/F}$
\item For $x \in C_K \subset W_{K/F}$, one has 
$$\alpha_\psi(x) = \psi(x/x^\sigma).$$
\item For $x = \widetilde{\sigma} y \in  W_{K/F} \setminus C_K$, one has 
$$\alpha_\psi(x) = \psi(y^\sigma/y) \rtimes \sigma.$$
\end{itemize}
To check this, one needs to make explicit calculations in group homology and cohomology, for which we use the usual ``inhomogeneous'' free resolution of $\Z$. Since the necessary formulas on the homology side are not as widely standardized as those on the cohomology side, we list them explicitly.

Let $G$ be a group and $M$ a left $G$-module. Computing via the inhomogeneous resolution gives the usual description
$$
H^1(G, M) = \frac{\{ \text{maps } \xi: G \to M \text{ satisfying } \xi(gh) = \xi(g) + g\xi(h) \}}{ \{ \text{maps such that there exists } m \in M, \xi(g) = gm - m\}}
$$
If $\oplus_S N$ is a direct sum of copies of an abelian group $N$ indexed by a set $S$, the symbol $\delta_s(n)$  refers to the element which is $n$ in the spot indexed by $s$ and $0$ elsewhere. Computing via the inhomogeneous resolution then gives
$$
H_1(G, M) = \frac{\{(m_g)_{g \in G} | \sum_g (g^{-1}m_g - m_g) = 0\}}{ d(\oplus_{G \times G} M)},
$$
where $d(\delta_{g, h}(m)) = \delta_h(g^{-1}m) - \delta_{gh}(m) + \delta_g(m)$. The cup product pairing
$$
H^1(G, N) \otimes H_1(G, M) \to M \otimes_{\mathbb{Z}[G]} N
$$
is given by $\xi \otimes (m_g) \mapsto \sum_g m_g \otimes \xi(g)$. If $G' < G$ is a finite-index normal subgroup, one has an action of $G/G'$ on $H_1(G', M)$ by the rule $g * \delta_{g'}(m) = \delta_{gg'g^{-1}}(gm)$, and the trace map 
$$
H_1(G, M) \to H_1(G', M)^{G/G'},
$$
may be computed as follows: pick coset representatives $g_1, g_2, \ldots, g_n$ for $G/G'$. Then any $g \in G$ determines a permutation $\tau \in S_n$ by the rule $g_i g = g' g_{\tau(i)}$ (where $g' \in G'$), and
$$
\text{Trace} (\delta_g(m)) = \sum_i \delta_{g_i g g_{\tau(i)}^{-1}}(g_i m).
$$

If $A$ is an abelian group and $G$ is an arbitrary group with a surjection to $\mathbb{Z}/2\Z$, write $\boxed{A}$ for the $G$-module with underlying group $A$ such that $G$ acts via $\mathbb{Z}/2\Z$ via inversion, and write $\boxed{a} \in \boxed{A}$ for the element $a \in A$. Then one has $L = \widehat{L} = \boxed{\Z}$, and $\widehat{T} = \boxed{\C^\times}$.

Now given a character $\psi$ of $T(F)^! = C_K^{(1)}$, to complete the calculation, we must show that 
$$
\xi_\psi(x) = \begin{cases} 
\psi(x/x^\sigma)& \text{ for } x \in C_K \\
\psi((x\widetilde{\sigma}^{-1})/(\widetilde{\sigma}^{-1} x)) & \text{ for } x \in W_{K/F} \setminus C_K.
\end{cases}
$$
is the cohomology class given by the recipe above. 
To prove this claim, given $\alpha \in T(F)^!$, pick $x \in C_K$ with $x/x^\sigma = \alpha$.  With the explicit formulas for trace and cup product above, one computes:
\begin{itemize}
\item that identifications (\ref{standard1})-(\ref{standard3}) send $\alpha$ to the 1-cycle $\delta_\alpha(\boxed{1})$.
\item that the image of $\delta_x(\boxed{1})  \in H_1(W_{K/F}, \widehat{L})$ under the trace $H_1(W_{K/F}, \widehat{L}) \to H_1(C_K, \widehat{L})^{\Gal(K/F)}$ is $\delta_\alpha(\boxed{1})$.
\item that $\xi_\psi$ satisfies the cocycle condition (using $\widetilde{\sigma}^2 \in C_F$). 
\item that the cup product $(\delta_x(\boxed{1}), \xi_\psi)$ is $\boxed{\psi(\alpha)} \otimes \boxed{1} = \psi(\alpha)$.
\end{itemize}

\subsection*{Acknowledgements}
This research was supported by the Swiss National Science Foundation. We thank Eva Bayer, Farrell Brumley, Kevin Buzzard, Emmanuel Kowalski, Dimitar Jetchev, Ting-Yu Lee, Philippe Michel, Djordje Mili\'cevi\'c, and Julian Rosen for helpful conversations.

\def\cprime{$'$}

\end{document}